\documentclass{amsart}
\usepackage{graphicx}
\usepackage{colordvi}
\newtheorem{thm}{Theorem}[section]
\newtheorem{cor}[thm]{Corollary}
\newtheorem{lem}[thm]{Lemma}
\newtheorem{prop}[thm]{Proposition}
\theoremstyle{definition}
\newtheorem{defn}[thm]{Definition}
\theoremstyle{remark}
\newtheorem{rem}[thm]{Remark}
\numberwithin{equation}{section}
\newcommand{\norm}[1]{\left\Vert#1\right\Vert}
\newcommand{\abs}[1]{\left\vert#1\right\vert}

\newcommand{\Real}{\mathbb R}

\newcommand{\pfrac}[2]{\frac{\partial #1}{\partial #2}}
\title[Ricci flow on conical surfaces]{Ricci flow on surfaces surfaces with conical singularities}

\author{Hao Yin}

\thanks{Mathematics Department, Shanghai Jiaotong University, 200240, Shanghai, China. E-mail: haoyin@sjtu.edu.cn.}

\begin{document}

\vspace{5cm}
{\centerline { \Huge Ricci flow on surfaces with conical singularities}}

\vspace{3cm}
{\centerline{\Huge Hao Yin}}

\vfill
\noindent
Running title:Ricci flow on conical surfaces\\
Address: Mathematics Department, Shanghai Jiaotong University, 200240, Shanghai, China.\\
Email: haoyin@sjtu.edu.cn\\

\newpage
\Large
\noindent
{\bf Abstract.}

This paper studies the normalized Ricci flow on surfaces with conical singularities. It's proved that the normalized Ricci flow has a solution for a short time for initial metrics with conical singularities. Moreover, the solution makes good geometric sense. For some simple surfaces of this kind, for example, the tear drop and the football, it's shown that they admit Ricci soliton metric.

\vspace{2cm}

\vspace{2cm}


\newpage
\section{Introduction}

Let $S$ be a closed Riemann surface. It's well known that there exists a collection of Riemannian metrics $g$ compatible with the conformal structure of $S$. This collection of metrics is called a conformal class of metrics. By Uniformization theorem, there exists a special metric in this class such that the Gaussian curvature is constant, whose sign is given by Euler number of $S$. 
Moreover, in 1982, Hamilton \cite{Ham} used his Ricci flow to prove that each metric in a given conformal class can be deformed into a metric of constant curvature. He assumed some technical condition if the Euler number is positive, which is later removed by Chow \cite{Ben}. However, their proofs need the Uniformization theorem in the case of positive Euler number. Very recently, a short paper of Chen, Lu and Tian \cite{CLT} made it clear that Ricci flow can be used to give a proof of the Uniformization of Riemann surfaces.

In this paper, we are concerned with metrics with some conical singularities. Precisely, let $S$ be a smooth Riemann surface and $p\in S$. Let $U$ be a neighborhood of $p\in S$ and $z$ be the conformal coordinate in $U$ such that $z(p)=0$. A metric $g$ on $S$ is said to have a \textit{conical singularity} of order $\beta(\beta>-1)$ , or of angle $\theta =2\pi(\beta+1)$ at $p$ if in $U$,
\begin{equation*}
g=\rho(z)\abs{z}^{2\beta}\abs{dz}^2,
\end{equation*}
for some continuous positive function $\rho$ in $U$. 

Now given $n$ points $p_i, i=1,\cdots,n$ on $S$ and $n$ numbers $\beta_i$ bigger than $-1$, there is the collection of singular metrics on $S$, each of which has a conical singularity of order $\beta_i$ at $p_i$ and is smooth and compatible with the conformal structure of $S$ elsewhere. A metric $g$ in this collection is called a conformal metric on $S$ with conical singularities prescribed by the divisor $\beta=\sum_{i=1}^n \beta_i p_i$. It turns out there is also an Euler number
\begin{equation*}
\tilde{\chi}=\chi(S)+\sum_{i=1}^n \beta_i
\end{equation*}
and a Gauss-Bonnet formula
\begin{equation}\label{Gauss_Bonnet}
\frac{1}{2\pi}\int_S K_g dA_g=\tilde{\chi},
\end{equation}
where $A_g$ and $K_g$ are the area element and the Gauss curvature of $g$. Please note that as given by Proposition 1 of \cite{Troyanov}, one needs to assume that $K_g$ as a function on $S$ is H\"{o}lder continuous in order that the Gauss-Bonnet formula is true. As we will see, in later discussion, this condition can be weakened. 

The first aim of this paper is to justify that there is a reasonable notion of normalized Ricci flow. That is, for some 'good' $g$ with prescribed conical singularities $\beta$, there is a constant $T>0$ and a smooth family $g(t), t\in [0,T]$ of conformal metrics with prescribed conical singularities $\beta$ such that $g(0)=g$ and
\begin{equation}\label{eqn_normalizedFlow}
\frac{\partial g}{\partial t}=(r-R)g.
\end{equation}
Here $R$ is the scalar curvature of $g(t)$ and $r$ is a constant. It will be shown that the Gauss-Bonnet formula is true for $t>0$ and if we choose $r$ to be the average of scalar curvature, then the volume of $g(t)$ is a constant.
For the precise formulation of the main result, we need a different (more regular) definition of conformal metric with conical singularities. Therefore, we will postpone the statement of the first result.

A natural question is that whether this flow exists for long time or not and if it does, whether it will converge to some special metric. Due to technical difficulties, the author is not able to answer this question yet. 

A closely related question is what is the most natural metric on a surface with conical singularities. The most natural 'special' metric is a metric with constant Gaussian curvature. It's well known that for some pairs of $(S,\beta)$, there doesn't exist a metric of constant curvature, while for some other pairs $(S,\beta)$, there does. There are some results on this question, but a necessary and sufficient condition is not known yet. See the discussion in \cite{Troyanov}, \cite{LuoTian} and \cite{threePoint}. 

When the metric of constant curvature doesn't exist, Chen \cite{ChenXX} proposed the extremal Hermitian metric, which is the critical point of Calabi functional, as a candidate of the 'natural' metric. A special case of the extremal Hermitian metric other than the constant curvature metric is the HCMU metric, meaning 'Hessian of curvature of the metric is umbilical'. Many authors studied the property of such metrics and constructed several examples. See for example, \cite{ChenXX}, \cite{LinZhu} and \cite{HCMU}. 

Now that since the normalized Ricci flow makes good sense, there is another candidate for the 'natural' metric if the constant curvature metric doesn't exist, i.e. the Ricci soliton. The other aim of this paper is to show that there exists, in a reasonable sense, a Ricci soliton metric for some simple surfaces with conical singularities known as the tear drop and the football.

The organization of the paper is as follows. In Section \ref{two}, we define some weighted H\"{o}lder space and pave the way for later analysis. In Section \ref{three}, we study the linear heat equation on conical surfaces. An existence theorem and some basic estimates are proved. In Section \ref{four}, Schauder fixed point theorem is applied to study the local existence of the Ricci flow equation. We also argue that the solution makes good geometric sense. In the last section, we use ODE method to show the existence of some Ricci soliton metric in the case of Riemann sphere with one or two conical singularities.

\textbf{Acknowledgment.} The author would like to thank Professor Chun-qin Zhou for introducing to him the topic of conical singularity and Professor De-liang Xu for numerous discussions on this paper.

\section{Weighted H\"{o}lder space}\label{two}

The pair $(S,\beta)$ of a smooth closed Riemann surface $S$ and a divisor $\beta=\sum \beta_i p_i$ contains the information of a singular conformal(smooth) structure. The purpose of this section is to make it clear what we mean by weighted H\"{o}lder continuous function with respect to this singular smooth structure in a coordinate-independent way. The space of weighted H\"{o}lder functions plays an important role in our later proof. The definition comes up naturally in the process of doing interior estimate for the Laplacian operator associated with cone metrics. 

Let's start from a local point of view. Assume that $p$ is a cone singularity of order $\beta$. In a neighborhood $U$ of $p$, there is local coordinates $(x,y)$, which is smooth with respect to $S$. Polar coordinates are given by
\begin{equation*}
	x=r\cos \theta,\quad y=r\sin \theta.	
\end{equation*}
Set 
\begin{equation*}
	\rho=\frac{1}{\beta+1}r^{\beta+1}.
\end{equation*}
\begin{rem}
	The reason for considering this transformation is that it transforms $r^{2\beta}(dr^2+r^2\theta^2)$ into $d\rho^2+\rho^2(\beta+1)^2 d\theta^2$. This will be explained more in the definition of weighted H\"{o}lder continuous metric on $(S,\beta)$.
\end{rem}

Given any function $f(x,y)=f(r,\theta)$, it is at the same time a function of $(\rho,\theta)$
\begin{equation*}
	F(\rho,\theta)=f( ( (\beta+1)\rho)^{\frac{1}{\beta+1}},\theta).
\end{equation*}
We will say that $f(x,y)$ is a weighted H\"{o}lder continuous function if $F(\rho,\theta)$ satisfies the assumptions below.

Assume without loss of generality that $F$ is defined on $0<\rho\leq 1$ and $\theta\in S^1$. For each $k\in \mathbb{N}$, set
\begin{equation*}
	\Omega_k=\{ (\rho,\theta)|\quad \frac{1}{2^{k+1}}<\rho<\frac{1}{2^{k-1}}\}.
\end{equation*}
Denote the restriction of $F$ to $\Omega_k$ by $F_k$ and consider a further substitution
\begin{equation*}
	s=2^{k} \rho.
\end{equation*}
$F_k$ is then a function of $(s,\theta)$ defined on $(\frac{1}{2},2)\times S^1$. Let $\norm{\cdot}_{C^{l,\alpha}}$ be the usual H\"{o}lder norm on this domain. 
\begin{defn}\label{defn_local}
	The $C^{l,\alpha}$ weighted H\"{o}lder continuous function $f(x,y)$ on $U$ is defined to be a function $f$ such that the following norm is finite,
	\begin{equation*}
		\norm{f}_{C^{l,\alpha}_b(U)}:= \sup_{k\in \mathbb{N}} \norm{F_k(s,\theta)}_{C^{l,\alpha}}.
	\end{equation*}
\end{defn}

We then move on to show that this definition is independent of the particular choice of smooth coordinates $(x,y)$ in $U$. For this purpose, let's give an equivalent definition of weighted H\"{o}lder continuous. Set
\begin{equation*}
w=\log_2 \rho
\end{equation*}
and
\begin{equation*}
G(w,\theta)=F(2^w,\theta).
\end{equation*}
\begin{prop}\label{defn_equiv}
$f$ is a weighted $C^{l,\alpha}$ function on $U$ if and only if $G(w,\theta)$ as a function on $(-\infty,0)\times S^1$ has uniform $C^{l,\alpha}$ norm on any tube $(-k-1,-k+1)\times S^1$.
\end{prop}
\begin{proof}
It suffices to show that 
\begin{equation*}
\frac{1}{C}\norm{F_k(s,\theta)}_{C^{l,\alpha}( (1/2,2)\times S^1)}<\norm{G(w,\theta)}_{C^{l,\alpha}( (-k-1,k+1)\times S^1)}< C \norm{F_k(s,\theta)}_{C^{l,\alpha}( (1/2,2)\times S^1)}
\end{equation*}
for some constant independent of $k$. This is true by observing that
\begin{equation*}
G(w,\theta)=F_k(2^{w+k},\theta)
\end{equation*}
and the function $w\mapsto 2^{w+k}$ defined on $(-k-1,-k+1)$ has bounded $C^l$ norm uniformly for all $k$.
\end{proof}

Consider a smooth coordinate transformation
\begin{equation*}
\left\{
\begin{array}[]{rcl}
\tilde{x}&=& \tilde{x}(x,y)\\
\tilde{y}&=& \tilde{y}(x,y)
\end{array}
\right.
\end{equation*}
such that $(0,0)$ is mapped to itself. It induces a transformation from $(w,\theta)$ to $(\tilde{w},\tilde{\theta})$. 

\begin{rem}
Here we pretend that $\theta$ and $\tilde{\theta}$ are functions while they take values in $S^1$. 
\end{rem}

\begin{eqnarray*}
\pfrac{}{w}&=& \pfrac{\rho}{w}\pfrac{}{\rho} \\
&=& (\log 2) \rho\pfrac{}{\rho}.
\end{eqnarray*}
Recall that
\begin{equation*}
\rho=\frac{1}{\beta+1}r^{\beta+1}.
\end{equation*}
Hence
\begin{eqnarray*}
\pfrac{}{w}&=& \frac{\log 2}{\beta+1}r^{\beta+1}\pfrac{r}{\rho}\pfrac{}{r} \\
&=& \frac{\log 2}{\beta+1}r\pfrac{}{r}\\
&=& \frac{\log 2}{\beta+1}(x\pfrac{}{x}+y\pfrac{}{y})\\
&=& \frac{\log 2}{\beta+1}\left( (x\pfrac{\tilde{x}}{x}+y\pfrac{\tilde{x}}{y})\pfrac{}{\tilde{x}} + (x\pfrac{\tilde{y}}{x}+y\pfrac{\tilde{y}}{y})\pfrac{}{\tilde{y}} \right).
\end{eqnarray*}
\begin{eqnarray*}
	\pfrac{}{\theta}&=&  \pfrac{x}{\theta}\pfrac{}{x}+\pfrac{y}{\theta}\pfrac{}{y} \\
	&=& -y\pfrac{}{x}+x\pfrac{}{y} \\
	&=& (x\pfrac{\tilde{x}}{y}-y\pfrac{\tilde{x}}{x})\pfrac{}{\tilde{x}}+ (x\pfrac{\tilde{y}}{y}-y\pfrac{\tilde{y}}{x})\pfrac{}{\tilde{y}})
\end{eqnarray*}
On the other hand,
\begin{equation*}
\pfrac{}{\tilde{w}}=\frac{\log 2}{\beta+1}(\tilde{x}\pfrac{}{\tilde{x}}+\tilde{y}\pfrac{}{\tilde{y}})
\end{equation*}
and
\begin{equation*}
\pfrac{}{\tilde{\theta}}=-\tilde{y}\pfrac{}{\tilde{x}}+\tilde{x}\pfrac{}{\tilde{y}}.
\end{equation*}
It's obvious that the constant $\frac{\log 2}{\beta+1}$ in the above equations plays no role in the following proof. For the sake of simplicity, we just pretend that there is no such constant. Next, we would like to write $\pfrac{}{w}$ and $\pfrac{}{\theta}$ as linear combinations of $\pfrac{}{\tilde{w}}$ and $\pfrac{}{\tilde{\theta}}$,
\begin{equation*}
	\left\{
	\begin{array}[]{rcl}
		\pfrac{}{w}&=& A\pfrac{}{\tilde{w}} + B\pfrac{}{\tilde{\theta}} \\
		\pfrac{}{\theta} &=&  C\pfrac{}{\tilde{w}} + D\pfrac{}{\tilde{\theta}}.
	\end{array}
	\right.
\end{equation*}
Therefore
\begin{equation*}
	\frac{\partial^m}{\partial w^m}\frac{\partial^n}{\partial \theta^n}=(A\pfrac{}{\tilde{w}} + B\pfrac{}{\tilde{\theta}})^m(C\pfrac{}{\tilde{w}} + D\pfrac{}{\tilde{\theta}})^n
\end{equation*}

To establish equivalence between $C^{l,\alpha}$ norm in $(w,\theta)$ and $C^{l,\alpha}$ norm in $(\tilde{w},\tilde{\theta})$, it suffices to prove 

$(\star)\quad $ $A,B,C,D$ together with their derivatives with respect to $\pfrac{}{\tilde{w}}$ and $\pfrac{}{\tilde{\theta}}$ are bounded uniformly.

By direct computation,
\begin{equation*}
	A=\frac{1}{\tilde{x}^2+\tilde{y}^2}(x\tilde{x}\pfrac{\tilde{x}}{x}+\tilde{x}y\pfrac{\tilde{x}}{y}+\tilde{y}x\pfrac{\tilde{y}}{x}+\tilde{y}y\pfrac{\tilde{y}}{y}).
\end{equation*}
The expressions for $B,C,D$ are similar. Each term in $A$ and its derivatives with respect to $\pfrac{}{\tilde{w}}$ and/or $\pfrac{}{\tilde{\theta}}$ has the following properties,
\begin{itemize}
	\item the denominator is a power of $\tilde{x}^2+\tilde{y}^2$,

	\item the numerator can be regarded formally as polynomial of $x,y,\tilde{x},\tilde{y}$ with coefficients given by multiplication of derivatives of $\tilde{x},\tilde{y}$ with respect to $x,y$ and derivatives of $x,y$ with respect to $\tilde{x},\tilde{y}$,

	\item the degree of numerator is no less than the degree of denominator as a polynomial of $\tilde{x}$ and $\tilde{y}$.
\end{itemize}
As an example, one can see such properties are true in the expansion of
\begin{equation*}
	\tilde{x}\pfrac{}{\tilde{x}}\left( \frac{\tilde{x}x\pfrac{\tilde{y}}{x}}{\tilde{x^2}+\tilde{y}^2} \right).
\end{equation*}
Similar properties are true for $B,C,D$. With this in mind, $(\star)$ follows from the following facts,
\begin{itemize}
	\item the derivatives of $\tilde{x},\tilde{y}$ with respect to $x,y$ are uniformly bounded (when $x,y$ are small) and the same is true for the derivatives of $x,y$ with respect to $\tilde{x},\tilde{y}$.

	\item There exists a constant $C$ such that $\frac{1}{C}(x^2+y^2)<\tilde{x}^2+\tilde{y}^2<C(x^2+y^2)$.

	\item $x<\sqrt{x^2+y^2}$ and so on.
\end{itemize}

In summary, we proved Definition \ref{defn_local} is independent of the choice of local coordinates. 

For each $p_i$, choose neighborhoods $U_i$ and coordinates $(x_i,y_i)$ as before so that $\norm{\cdot}_{C_b^{l,\alpha}(U_i)}$ makes sense. Away from these singular points, choose regular coordinate charts $V_j$ such that $U_i$ and $V_j$ cover $S$.
\begin{defn}
\label{defn_global}
Given $U_i$ and $V_j$ as above, define $C^{l,\alpha}(S,\beta)$ to be the set of functions with finite weighted $C^{l,\alpha}$ norm defined as following,
\begin{equation*}
\norm{f}_{C^{l,\alpha}(S,\beta)}=\max_i \norm{f}_{C_b^{l,\alpha}(U_i)} +\max_j \norm{f}_{C^{l,\alpha}(V_j)}.
\end{equation*}
\end{defn}

If $S$ is a closed 2-dimensional Riemannian manifold, then the pair $(S,\beta)$ represents a 'singular' smooth structure so that $C^{l,\alpha}(S,\beta)$ makes sense. If $S$ is a closed Riemann surface, then we have the additional conformal structure. We will now make it clear what is $C^{l,\alpha}$ cone metric in this conformal class $(S,\beta)$.

\begin{defn}
	\label{defn_metric}
	For each  $p\in S$, let $(x,y)$ be conformal coordinates around $p$ of $S$. Let $g$ be a metric on $S$. $g$ is said to be a $C^{l,\alpha}$ cone metric in the conformal class $(S,\beta)$ if for each nonsingular point $p$, $g=e^{2u}(dx^2+dy^2)$ for some $C^{l,\alpha}$ $u$ in usual sense and for each singular $p=p_i$, $g=e^{2u}(x^2+y^2)^{\beta_i}(dx^2+dy^2)$ for some weighted $C^{l,\alpha}$ $u$ as given in Definition \ref{defn_local}. 
	\end{defn}

If $p$ is a singular point and $(u,v)$ is another conformal coordinate system, then
\begin{equation*}
	(u^2+v^2)^\beta(du^2+dv^2)=\frac{(u^2+v^2)^\beta}{(x^2+y^2)^\beta}(\pfrac{u}{x}\pfrac{v}{y}-\pfrac{u}{y}\pfrac{v}{x})(x^2+y^2)^\beta(dx^2+dy^2).
\end{equation*}
Moreover, $\frac{(u^2+v^2)^\beta}{(x^2+y^2)^\beta}(\pfrac{u}{x}\pfrac{v}{y}-\pfrac{u}{y}\pfrac{v}{x})$ is a weighted $C^{l,\alpha}$ function. Therefore, Definition \ref{defn_metric} is independent of the choice of coordinates.

We also need some parabolic counterpart of the above weighted H\"{o}lder functions, since we will be working on parabolic equations. Let $p$ be a cone singularity of order $\beta$. Let $(x,y)$ be conformal coordinates in a neighborhood $U$ of $p$ and $(r,\theta)$ be the polar coordinates. For a function $f:U\times [0,T]\to \Real$, one may define as before
\begin{equation*}
	F(\rho,\theta,t)=f( ( (\beta+1)\rho)^{\frac{1}{\beta+1}},\theta,t).
\end{equation*}

Define $\Omega_k$ as before. Denote the restriction of $F$ to $\Omega_k\times[0,T]$ by $F_k$ and consider a further substitution
\begin{equation*}
	s=2^k\rho \mbox{ and } \tilde{t}=2^{2k} t.
\end{equation*}
$F_k$ is then a function of $(s,\theta,\tilde{t})$ defined on $(1/2,2)\times S^1\times [0,2^{2k}T]$. Let $\norm{\cdot}_{C^{l,\alpha}(\Omega\times [0,T])}$ be the space-time H\"{o}lder space defined in Chapter IV Section 1 of \cite{Lie}. Now we can define
\begin{defn}\label{defn_paralocal}
	The $C^{l,\alpha}$ weighted H\"{o}lder continuous function $f(x,y,t)$ on $U\times [0,T]$ is defined to be a function $f$ such that the following norm is finite
	\begin{equation*}
		\norm{f}_{C^{l,\alpha}_{b}(U\times [0,T])}:=\sup_{k\in \mathbb{N}}\norm{F_k(s,\theta,\tilde{t})}_{C^{l,\alpha}( (1/2,2)\times S^1\times [0,2^{2k}T])}.
	\end{equation*}
\end{defn}

Choose $U_i$ and $V_i$ as before. We can now define space-time H\"{o}lder space globally.
\begin{defn}
\label{defn_paraglobal}
Define $C^{l,\alpha}( (S,\beta)\times [0,T])$ to be the set of functions with finite weighted $C^{l,\alpha}$ norm defined as following,
\begin{equation*}
\norm{f}_{C^{l,\alpha}( (S,\beta)\times [0,T])}=\max_i \norm{f}_{C_b^{l,\alpha}(U_i\times [0,T])} +\max_j \norm{f}_{C^{l,\alpha}(V_j\times [0,T])}.
\end{equation*}
\end{defn}

\begin{rem}
	Due to almost the same proof, this definition of $C^{l,\alpha}( (S,\beta)\times [0,T])$ doesn't depend on the choice of coordinate systems.
\end{rem}

\section{Linear equations}\label{three}

Let $g$ be a $C^{l,\alpha}$ cone metric in the conformal class $(S,\beta)$. In this section, we are concerned with the linear parabolic equation
\begin{equation}
	\pfrac{u}{t}=a(x,t)\triangle u(x,t)+f(x,t) 
	\label{eqn_linear}
\end{equation}
where $\triangle$ is the Laplacian operator of metric $g$. We will prove that given $a$ and $f$ in $C^{l,\alpha}( (S,\beta)\times [0,T])$, for each $u_0$ in $C^{l,\alpha}( (S,\beta))$, there exists a solution $u(x,t)$ to equation (\ref{eqn_linear}) with initial value $u_0$. Moreover, the solution satisfies some estimates.

Before we get into the exact statement of these results, let's recall that since the metric is not smooth everywhere, we start by defining the Laplacian on smooth functions supported away from the singularities and then extend the domain of definition to get an self-adjoint operator. It's well known that there are many different self-adjoint extensions in the case of surfaces with conical singularities. Moreover, different self-adjoint extension means different parabolic equation. There is an implicit choice of boundary data at the singular points. Among these choices, there is the one which makes our geometric problem meaningful. The choice will be implicit in the construction of solution given below. We will show later that it's compatible with our geometric problem.

\begin{thm}\label{thm_linear}
	For each $a,f$ in $C^{l,\alpha}( (S,\beta)\times [0,T])$ and $u_0$ in $C^{l+2,\alpha}(S,\beta)$, there exists a solution $u(x,t)$ to equation (\ref{eqn_linear}) with initial value $u(x,0)=u_0(x)$. Moreover,
	\begin{equation*}
		\norm{u(x,t)}_{C^{l+2,\alpha}( (S,\beta)\times [0,T])}\leq C(\norm{a}_{C^{l,\alpha}( (S,\beta)\times [0,T])}, \norm{f}_{C^{l,\alpha}( (S,\beta)\times [0,T])},\norm{u_0}_{C^{l+2,\alpha}( S,\beta)}).
	\end{equation*}
\end{thm}

	Assume without loss of generality that there is only one conical singular point $p$ of order $\beta$. Fix a conformal coordinate $(x,y)$ in a neighborhood $U$ of $p$. Let $(r,\theta)$ be the polar coordinates and set
	\begin{equation*}
		\rho=\frac{1}{\beta+1}r^{\beta+1}
	\end{equation*}
	as before. Consider surfaces with boundary
	\begin{equation*}
		S_k=S\setminus\{(x,y)\in U| \rho(x,y)<2^{-k}\}.	
	\end{equation*}
	Since $a,f$ and $u_0$ are all $C^{l,\alpha}$ on $S_k$ up to the boundary, it's well known that there exists a solution $u_k(x,t)$ defined on $S_k\times [0,T]$ to the following equation
	\begin{equation}
		\left\{
		\begin{array}{ll}
			\pfrac{u_k}{t}=a(x,t)\triangle u_k+f(x,t) & \mbox{in} S_k \\
			\pfrac{u_k}{\nu}=0 & \mbox{on} \partial S_k
		\end{array}
		\right.
		\label{eqn_sk}
	\end{equation}
	where $\nu$ is the outward normal to the boundary.

	We will prove that when $k$ goes to infinity, $u_k$ will converge to a solution of equation (\ref{eqn_linear}). For this purpose, we need the following maximum principle to provide $C^0$ a priori estimate.

\begin{prop}\label{prop_maximum} 
	Assume that $a(x,t)$ is nonnegative. $\abs{f(x,t)}$ is uniformly bounded by $C_2$. $\abs{u_0}$ is bounded by $C_1$. Then for $t\in [0,T)$,
	\begin{equation*}
		\max \abs{u_k(\cdot,t)}\leq C_1+C_2 t.
	\end{equation*}
\end{prop}
\begin{proof}
Consider first the upper bound. Given any small $\varepsilon>0$. Fix some smooth function $w$ defined on $S_k$ such that
\begin{equation*}
	\pfrac{w}{\nu}<0
\end{equation*}
on $\partial S_k$. Let's consider the equation satisfied by $\tilde{u}(x,t)=u_k(x,t)+\varepsilon w(x)$.
\begin{equation*}
	\pfrac{\tilde{u}}{t}=\pfrac{u_k}{t}=a(x,t)\triangle \tilde{u}+f(x,t)-\varepsilon a(x,t)\triangle w.
\end{equation*}
Set $\tilde{f}(x,t)=f(x,t)-\varepsilon a(x,t)\triangle w(x)$.
We have
\begin{equation*}
	\pfrac{\tilde{u}}{t}=a(x,t)\triangle \tilde{u}+\tilde{f}.
\end{equation*}
And,
\begin{equation*}
	\pfrac{\tilde{u}}{\nu}(x,t)<0,
\end{equation*}
for all $x\in \partial S_k$ and $t\in [0,T]$.
Hence,
\begin{equation*}
	\pfrac{\tilde{u}}{t}\leq a(x,t)\triangle \tilde{u} +C_2(\varepsilon).
\end{equation*}
Here $C_2(\varepsilon)$ depends on the choice of $w$. Moreover, when $\varepsilon$ goes to zero, $C_2(\varepsilon)$ converges to $C_2$.

Set $h(t)=\max_{x\in S_k} \tilde{u}(x,t)$. Due to the lemma in Section 3.5 of \cite{Ham4}, 
\begin{equation*}
	\frac{d}{dt}h(t)\leq \sup_{x\in Y(t)}\pfrac{\tilde{u}(y,t)}{t},
\end{equation*}
where $Y(t)$ is the set of all points $y$ such that $\tilde{u}(y,t)=\max_{x\in S_k}\tilde{u}(x,t)$. Due to the boundary condition of $\tilde{u}$, $Y(t)$ is contained in the interior of $S_k$. Hence for each $y\in Y(t)$,
\begin{equation*}
	(\triangle \tilde{u})(y,t)\leq 0.
\end{equation*}
Therefore 
\begin{equation*}
	\frac{dh}{dt}\leq C_2(\varepsilon),
\end{equation*}
from which we have
\begin{equation*}
	\tilde{u}(x,t)\leq C_1+\varepsilon \max_{x\in S_k}\abs{w(x)} +t C_2(\varepsilon).
\end{equation*}

Let $\varepsilon$ go to zero.
\begin{equation*}
	u_k(x,t)\leq C_1+C_2 t.
\end{equation*}
The proof for lower bound is similar.
\end{proof}
\begin{rem}
	The above result doesn't depend on the geometry of $\partial S_k$. Once $f$ and $u_0$ are uniformly bounded on $S$ as assumed in Theorem \ref{thm_linear}, we have uniform $C^0$ estimate for all $u_k$.
\end{rem}

We now move on to show higher order estimates of $u_k$. As one may expect, due to the existence of the singular point, the closer we are to the singularity, the worse the estimate is. That explains partially why we use weighted H\"{o}lder space. The main tool here is the H\"{o}lder interior estimates in standard theory of parabolic equations.  Our main reference will be the book of Lieberman \cite{Lie}. For the reader's convenience, we will give the exact statement below. 
\begin{prop}\label{prop_holder}
Let $\Omega$ be some domain in $\Real^n$ and $u(x,t)$ be a solution to
\begin{equation*}
	\pfrac{u}{t}=\sum_{i,j} a_{ij}(x,t)\partial_{ij}u +\sum_i b_i(x,t)\partial_i u+ f(x,t)
\end{equation*}
on $\Omega\times [0,T]$ with initial value $u(x,0)=u_0(x)$. Assume that
\begin{equation*}
	\lambda\abs{\xi}^2\leq a_{ij}\xi^i\xi^j\leq \Lambda\abs{\xi}^2.
\end{equation*}
Then for any compact set $K\subset \Omega$, we have
\begin{equation*}
	\norm{u}_{C^{l+2,\alpha}(K\times[0,T])}\leq C(\norm{u}_{C^0(\Omega\times [0,T])}+\norm{f}_{C^{l,\alpha}(\Omega\times [0,T])}+\norm{u_0}_{C^{l+2,\alpha}(\Omega)})
\end{equation*}
where $C$ depends on 
$\Omega,K,\lambda,\Lambda,\norm{a_{ij}}_{C^{l,\alpha}(\Omega\times [0,T])},\norm{b_j}_{C^{l,\alpha}(\Omega\times [0,T])}$.
\end{prop}

It's now time to prove Theorem \ref{thm_linear}. Recall that $S$ is covered by coordinate neighborhoods $U$ and $V_i$. Due to Proposition \ref{prop_maximum} and Proposition \ref{prop_holder}, we have
\begin{equation}\label{eqn_vi}
	\norm{u_k}_{C^{l+2,\alpha}(V_i\times [0,T])}\leq C(\norm{a}_{C^{l,\alpha}(V_i\times [0,T])},\norm{f}_{C^{l,\alpha}(V_i\times [0,T])},\norm{u_0}_{C^{l+2,\alpha}(V_i)}).
\end{equation}
To prove estimates in $U$, recall that $u_k(x,y)$, or $u_k(r,\theta)$ in polar coordinates may be regarded as a function of $\rho=\frac{1}{\beta+1}r^{\beta+1}$ and $\theta$. For the sake of simplicity, we will write $u_k(\rho,\theta)$ for it and this convention will be applied to other functions in the following proof. Moreover, in terms of $\rho$ and $\theta$, the cone metric will be
\begin{equation*}
	g=e^{2w}(d\rho^2+(\beta+1)^2\rho^2 d\beta^2),
\end{equation*}
where by Definition \ref{defn_metric} $w$ is in $C^{l,\alpha}_b(U)$.
Therefore, equation (\ref{eqn_linear}) becomes
\begin{equation*}
	\pfrac{u_k}{t}=a(\rho,\theta,t) e^{-2w(\rho,\theta)}(\frac{\partial^2}{\partial \rho^2}+\frac{1}{\rho}\pfrac{}{\rho}+\frac{1}{(\beta+1)^2\rho^2}\frac{\partial^2}{\partial \theta^2})u_k +f(\rho,\theta,t).
\end{equation*}

Fix some $m\in \mathbb{N}$. When $k$ is large so that $\Omega_m\subset S_k$, consider the restriction of $u_k$ to $\Omega_m\times [0,T]$ and a further substitution as given in Definition \ref{defn_paralocal},
\begin{equation*}
	s=2^{m}\rho,\mbox{ and } \tilde{t}=2^{2m}t.
\end{equation*}
For simplicity, we continue to write $u_k(s,\theta,\tilde{t})$ for this new function defined on $(1/2,2)\times S^1\times [0,2^{2m}T]$. The equation (\ref{eqn_linear}) becomes
\begin{equation*}
	\pfrac{u_k}{\tilde{t}}=a(s,\theta,\tilde{t})e^{-2w(s,\theta)}(\frac{\partial^2}{\partial s^2}+\frac{1}{s}\pfrac{}{s} +\frac{1}{(\beta+1)^2}\frac{\partial^2}{\partial\theta^2})u_k +2^{-2m}f(s,\theta,\tilde{t}).
\end{equation*}
By definition, $C^{l,\alpha}( (S,\beta)\times [0,T])$ norms of $a$ and $f$ dominate $C^{l,\alpha}$ norms of $a(s,\theta,\tilde{t})$ and $f(s,\theta,\tilde{t})$ on $(1/2,2)\times S^1\times [0,2^{2m}T]$. Recall that we have uniform $C^0$ estimate for $u_k$ so that we can now use Proposition \ref{prop_holder} to conclude
\begin{equation}\label{eqn_U}
	\norm{u_k}_{C^{l+2,\alpha}_b(U\times [0,T])}\leq C(\norm{a}_{C^{l,\alpha}_b(U\times [0,T])},\norm{f}_{C^{l,\alpha}_b(U\times [0,T])},\norm{w}_{C^{l,\alpha}_b(U\times [0,T])},\norm{u_0}_{C^{l+2,\alpha}_b(U)}).
\end{equation}

Due to equation (\ref{eqn_vi}) and (\ref{eqn_U}), a subsequence of $u_k$ converges to a solution $u(x,t)$ to equation (\ref{eqn_linear}). Moreover, we have
\begin{equation*}
	\norm{u}_{C^{l+2,\alpha}( (S,\beta)\times [0,T])}\leq C(\norm{a}_{C^{l,\alpha}( (S,\beta)\times [0,T])}+\norm{f}_{C^{l,\alpha}( (S,\beta)\times [0,T])}+\norm{u_0}_{C^{l+2,\alpha}( S,\beta)}).
\end{equation*}
This concludes the proof of Theorem \ref{thm_linear}.

Since the solution constructed in the proof of Theorem \ref{thm_linear} is the limit of $u_k$, the result of Proposition \ref{prop_maximum} holds for $u$. That is
\begin{equation*}
	\max \abs{u(\cdot,t)}\leq C_1+C_2 t
\end{equation*}
for all $t\in [0,T)$, where $C_1$ and $C_2$ bound $\abs{u_0}$ and $\abs{f}$ respectively. In rest of this paper, we will only use a special case of this result.

\begin{cor}
	\label{cor_special}
	If $u_0\equiv 0$ and $u(x,t)$ is the solution given by Theorem \ref{thm_linear}, then
	\begin{equation*}
		\abs{u(x,t)}\leq Ct,
	\end{equation*}
	where $C$ depends only on the $C^0$ norm of $f(x,t)$.
\end{cor}

The above results are enough to prove the local existence of the normalized Ricci flow. However, we need additional information on the solution $u$ to show that the solution of the normalized Ricci flow constructed in the next section makes good geometric sense. For example, one may expect that the volume is constant along the flow and that the Gauss-Bonnet theorem holds as long as the solution exists.

Let $u_k$ be the solution to equation (\ref{eqn_sk}) in the proof of Theorem \ref{thm_linear}.
\begin{eqnarray*}
	\frac{d}{dt}\int_{S_k} \abs{\nabla u_k}^2 dV 
	&=& 2 \int_{S_k} \nabla u_k\cdot \nabla(a(x,t)\triangle u_k +f (x,t)) dV\\
	&=& -2\int_{S_k} a(x,t)(\triangle u_k)^2 dV +2\int_{S_k} \nabla u_k\cdot \nabla f dV \\
	&\leq& \int_{S_k} \abs{\nabla u_k}^2 dV + \int_{S_k} \abs{\nabla f}^2 dV \\
	&\leq& \int_{S_k} \abs{\nabla u_k}^2 dV +C,
\end{eqnarray*}
where the integration by parts is justified by the boundary condition of $u_k$ and 
\begin{equation*}
	C=\max_{t\in [0,T]} \int_{S}\abs{\nabla f}dV.
\end{equation*}
ODE comparison gives
\begin{equation*}
	\int_{S_k}\abs{\nabla u_k(\cdot,t)}^2 dV \leq (\int_{S_k}\abs{\nabla u_0}^2 dV)e^t +C(e^t-1).
\end{equation*}
Let $k$ go to infinity. We have
\begin{lem}
	\label{lem_energy}
	Let $u(x,t)$ be the solution given by Theorem \ref{thm_linear} with $u_0\equiv 0$. Then
	\begin{equation*}
	\int_{S}\abs{\nabla u(\cdot,t)}^2 dV \leq (e^t-1) \max_{t\in [0,T]}\int_{S}\abs{\nabla f}^2 dV.
	\end{equation*}
\end{lem}

\section{Local existence}\label{four}

Let $g_0$ be a \textit{good} cone metric in the conformal class $(S,\beta)$. The main result of this section is to show there exists some $T_0>0$ such that the normalized Ricci flow has a solution defined on $[0,T_0]$ with initial metric $g_0$. Set $g(t)=e^{2u(x,t)}g_0$. The normalized Ricci flow equation becomes
\begin{equation}
	\pfrac{u}{t}=e^{-2u}\triangle u+ \frac{r}{2}-e^{-2u}K_0.
	\label{eqn_main}
\end{equation}
Here $K_0$ is the Gaussian curvature of $g_0$. For now, $r$ is just some constant. Later we will show that this $r$ can be taken as the average of the scalar curvature.

To make it clear that what we mean by a 'good' initial metric $g_0$, define a Banach space 
\begin{equation*}
W=C^{0,\alpha}(S,\beta)\cap H^1(S,\beta).
\end{equation*}
For $H^1(S,\beta)$, we need to fix some background $C^{0,\alpha}$ cone metric. It is not difficult to see that $H^1(S,\beta)$ is independent of the choice of this background metric. Throughout this section, we will assume that $g_0$ satisfies

(1) $g_0$ is a $C^{0,\alpha}$ cone metric in the conformal class $(S,\beta)$ as defined in Definition \ref{defn_metric}.

(2) The Gauss curvature $K_0$ of $g_0$ is in $W$.

The idea of proving local existence is the Schauder fixed point theorem, which is very common in the theory of nonlinear parabolic equations. We will use the following version of Schauder fixed point theorem,
\begin{thm}
	\label{thm_fixed}
	Let $E$ be a closed, convex set in a Banach space $V$, and let $F:E\to E$ be a continuous map such that $F(E)$ is relatively compact. Then $F$ has a fixed point.
\end{thm}

For some $T$ to be determined later, set
\begin{equation*}
V=C^{0,\alpha}( (S,\beta)\times [0,T])\cap C^0([0,T],H^1(S,\beta)).
\end{equation*}

\begin{rem}
The definition may look a little complicated. The idea is simple. That is we want our $u(x,t)$ to be in $C^{0,\alpha}( (S,\beta)\times [0,T])$ and to have finite energy $\int_M \abs{\nabla u}^2 dV$ for each $t\in [0,T]$.
\end{rem}

Set
\begin{equation*}
	E=\{u(x,t) \in V| u(x,0)\equiv0,  \norm{u}_{C^{0,\alpha}( (S,\beta)\times [0,T])}+\max_{t\in [0,T]} \norm{\nabla _{g_0} u}_{L^2(S,g_0)}\leq C_0 \}.
\end{equation*}
It is easy to check that $E$ is a closed and convex set in $V$. For $v\in E$, consider equation
\begin{equation}
	\pfrac{u}{t}=e^{-2v}\triangle u+\frac{r}{2}-e^{-2v}K_0.
	\label{eqn_schauder}
\end{equation}
Here $\triangle$ is the Laplacian operator of $g_0$. 

If we set $a(x,t)=e^{-2v}$ and $f(x,t)=\frac{r}{2}-e^{-2v}K_0$, then the $C^{0,\alpha}( (S,\beta)\times [0,T])$ norms of $a,f$ are bounded by the corresponding norms of $v$ and $K_0$. Theorem \ref{thm_linear} implies the existence of $u(x,t)$ which is a solution to equation (\ref{eqn_schauder}) with initial condition $u_0\equiv 0$. Moreover, the $C^{2,\alpha} ( (S,\beta)\times [0,T])$ of $u$ is bounded by $C^{0,\alpha}( (S,\beta)\times [0,T])$ norm of $v$. Hence, if we define a map $\Psi:E\to V$ by
\begin{equation*}
	\Psi(v)=u,
\end{equation*}
then $\Psi$ is a compact map. Now to prove the local existence of equation (\ref{eqn_main}), it suffices to show that when $T$ is small enough, $\Psi(v)\in E$. The fixed point given by Schauder fixed point theorem is the required solution.

\begin{lem}\label{lem_solution}
	There exists some $T>0$ such that for any $v\in E$, $\Psi(v)\in E$.
\end{lem}
\begin{proof}
	Let $u=\Psi(v)$. Due to Theorem \ref{thm_linear},
	\begin{equation}\label{eqn_high}
		\norm{u}_{C^{2,\alpha}( (S,\beta)\times [0,T])}\leq C_1 (C_0),
	\end{equation}
	where $C_1(C_0)$ is some constant depending on $C_0$. Due to Corollary \ref{cor_special},
	\begin{equation}\label{eqn_low}
		\norm{u}_{C^0 (S\times [0,T])}\leq C_2(C_0)T.
	\end{equation}
	Here $C_2(C_0)$ is another constant depending on $C_0$.

	We will show that if $T$ is small then 
	\begin{equation}\label{eqn_inter}
		\norm{u}_{C^\alpha ( (S,\beta)\times [0,T])}\leq \frac{1}{2} C_0.
	\end{equation}
	For this purpose, we need some interpolation inequality. Usually, a proof of interpolation inequality for H\"{o}lder space is very complicated. Fortunately, we need to estimate only the $C^\alpha$ norm for some $\alpha\in (0,1)$.  For the reader's convenience, we will show the proof below. Let $\Omega$ be any domain. Let $X=(x,t)$ and $Y=(y,s)$ be two space-time points in $\Omega\times [0,T]$. Let $\abs{X-Y}=\max \{\abs{x-y},\sqrt{\abs{t-s}}\}$.
	\begin{eqnarray*}
		\frac{\abs{u(X)-u(Y)}}{\abs{X-Y}^\alpha} &\leq& (2\norm{u}_{C^0(\Omega)})^{1-\alpha}\left( \frac{\abs{u(X)-u(Y)}}{\abs{X-Y}} \right)^\alpha .
	\end{eqnarray*}
	Hence,
	\begin{equation*}
		\norm{u}_{C^\alpha(\Omega \times [0,T])}\leq C \norm{u}_{C^0(\Omega\times [0,T])}^{1-\alpha} \norm{u}^\alpha_{C^{0,1}(\Omega \times [0,T])}.
	\end{equation*}
	It's important that this inequality holds uniformly regardless of the shape of $\Omega$ and $T$. We will then apply it to estimate $C^\alpha$ norm of u on $V_i\times [0,T]$ and $(1/2,2)\times S^1\times [0,2^{2m}T]$.  By the definition of $C^{l,\alpha}( (S,\beta)\times [0,T])$, one can choose $T$ so small that equation (\ref{eqn_inter}) is true.

	To complete the proof of the lemma, it remains to prove
	\begin{equation*}
		\max_{t\in [0,T} \norm{\nabla u}_{L^2(S,g_0)}\leq \frac{C_2}{2},
	\end{equation*}
	which is a direct consequence of Lemma \ref{lem_energy}.
\end{proof}

Now, Schauder fixed point theorem implies that there exists some $T>0$ and $u\in E$ such that $u$ is a solution to equation (\ref{eqn_main}) for $t\in [0,T]$. The next lemma improves the regularity of the solution.

\begin{lem}\label{lem_regularity}
	Let $g$ be $C^{l,\alpha}$ cone metric in the conformal class $(S,\beta)$. Let $u$ be the solution to equation (\ref{eqn_main}) obtained by Schauder fixed point theorem. Then there exists some constant $C$ such that
	\begin{equation*}
		\norm{u}_{C^{l+2}( (S,\beta)\times [0,T])}\leq C.
	\end{equation*}
\end{lem}
\begin{proof}
	The proof is repeated use of Proposition \ref{prop_holder} to $u$ on $V_i\times [0,T]$ and $(1/2,2)\times S^1\times [0,2^{2m}T]$.
\end{proof}

It now remains to show that this solution makes good geometric sense. That is, we want to show if we choose $r$ to be the average of scalar curvature for the initial metric, then the volume of $g(t)$ is constant and Gauss-Bonnet theorem still holds for $g(t)$. As a consequence, $r$ is the average of scalar curvature for $t>0$. We need this lemma,

\begin{lem}
	\label{lem_integration}
	For each $u\in E$,
	\begin{equation*}
		\int_S \triangle_{g(t)}u(x,t) dV_{g(t)}=0,
	\end{equation*}
	where $g(t)=e^{2u}g_0$.
\end{lem}
\begin{proof}
	Assume without loss of generality that there is only one cone singularity $p$. Let $(x,y)$ be a conformal coordinate system near $p$. Let $r,\theta$ and $\rho$ be as before. By integration by parts, it suffices to show
	\begin{equation*}
		\lim_{k\to \infty} \int_{\partial S_k}\abs{\nabla_{g(t)}u(x,t)_{g(t)}}d V_{g(t)}=0.
	\end{equation*}
	Let $g_0=e^{2w(x)}(d\rho^2+(\beta+1)^2 \rho^2 d\theta^2)$. Then
	\begin{equation*}
		g(t)=e^{2u(x,t)+2w(x)}(d\rho^2 + (\beta+1)^2\rho^2 d\theta^2).
	\end{equation*}
	Since both $u$ and $w$ are bounded, the length of $\partial S_k$ is less than $C 2^{-k}$. On the other hand, by the conformal invariance of energy and the fact that $u\in E$,
	\begin{equation*}
		\int_0^\varepsilon \int_0^{2\pi}\left( \abs{\pfrac{u}{\rho}}^2 +\frac{1}{\rho^2}\abs{\pfrac{u}{\theta}}^2 \right) \rho d\rho d\theta\leq C.
	\end{equation*}
	Let $s=\log_2 \rho$. The above equation is equivalent to
	\begin{equation*}
		\int_{-\infty}^{\log_2 \varepsilon} \int_0^{2\pi} \left( \abs{\pfrac{u}{s}}^2 +\abs{\pfrac{u}{\theta}}^2 \right) ds d\theta<C.
	\end{equation*}
	Hence,
	\begin{equation*}
		\lim_{k\to \infty}\int_{-k-1}^{-k+1} \int_0^{2\pi} \left( \abs{\pfrac{u}{s}}^2 +\abs{\pfrac{u}{\theta}}^2 \right) ds d\theta=0.
	\end{equation*}
	By Proposition \ref{defn_equiv}, $u$ as a function of $(s,\theta)$ is $C^{1,\alpha}$ in $(-k-1,-k+1)\times S^1$, which implies that
	\begin{equation*}
		\lim_{s\to -\infty} \abs{\pfrac{u}{s}}+\abs{\pfrac{u}{\theta}}= 0.
	\end{equation*}
	In terms of $(\rho,\theta)$, this is
	\begin{equation*}
		\lim_{\rho \to 0} \abs{\rho\pfrac{u}{\rho}}+\abs{\pfrac{u}{\theta}}=0.
	\end{equation*}
	Due to the boundedness of $u$ and $w$ again,
	\begin{equation*}
		\abs{\nabla_{g(t} u}_{g(t)}\leq C\sqrt{\abs{\pfrac{u}{\rho}}^2 +\frac{1}{\rho^2}\abs{\pfrac{u}{\theta}}^2}.
	\end{equation*}
	We can now estimate
	\begin{equation*}
		\int_{\partial S_k}\abs{\nabla_{g(t) u(x,t)}}dV_{g(t)}\leq C (\abs{\rho \pfrac{u}{\rho}}+\abs{\pfrac{u}{\theta}})|_{\rho=2^{-k}}\to 0.
	\end{equation*}

\end{proof}

Now we can prove the main result of this section.
\begin{thm}
	\label{thm_main} Let $g_0$ be a $C^{l,\alpha}$ cone metric in the conformal class of $(S,\beta)$ such that the Gauss curvature $K_0$ is in $W$. Assume that the Gauss-Bonnet theorem \ref{Gauss_Bonnet} is true for $g_0$. Choose $r$ to be the average of scalar curvature of $g_0$. Then there exists some $T>0$ and a solution $g(t)=e^{2u(x,t)}g_0 (0\leq t\leq T)$ to the normalized Ricci flow. Moreover, Gauss-Bonnet theorem is still true for $t>0$ and the volume of $g(t)$ is constant.
\end{thm}

\begin{proof}
By Theorem \ref{thm_fixed}, Lemma \ref{lem_solution} and Lemma \ref{lem_regularity}, we proved that there exists $T>0$ such that there is a $u(x,t)$ defined for $t\in [0,T]$, which is a solution of equation (\ref{eqn_main}). Obviously, $g(t)=e^{2u}g_0$ is a solution of the normalized Ricci flow (\ref{eqn_normalizedFlow}). It remains to prove that Gauss-Bonnet theorem is true for $t>0$ and the volume of $g(t)$ is constant. Although these are trivial facts for the normalized Ricci flow on smooth closed manifolds, they require some work in the case of singular surfaces.

Since Gauss-Bonnet is true for $g_0$, we have
\begin{equation}\label{aa}
\int_S K_0 dV_{g_0}=2\pi \tilde{\chi}.
\end{equation}
The equation for Gauss curvature under conformal change gives
\begin{equation}
e^{2u(x,t)}K_{t}=-\triangle_{g_0}u(x,t)+K_0.
\label{bb}
\end{equation}
Integrate equation (\ref{bb}) over $S$ with respect to the metric $g_0$,
\begin{equation*}
\int_S K_{t}dV_{g(t)}=-\int_{S}\triangle_{g_0}u(x,t)dV_{g_0} +\int_S K_0 dV_{g_0}.
\end{equation*}
It follows from equation (\ref{aa}) and Lemma \ref{lem_integration} that the Gauss-Bonnet is true for $g(t)$.

Let $V(t)$ be the volume of $g(t)$. Since we have chosen $r$ to be the average of scalar curvature of $g_0$,
\begin{equation}\label{cc}
V(0)=\frac{2}{r}\int_S K_0 dV_{g_0}.
\end{equation}
\begin{eqnarray*}
\frac{d}{dt}(V(t)-V(0)) &=& \frac{d}{dt}\int_S e^{2u}dV_{g_0} \\
&=&  2\int_S e^{2u}u_t dV_{g_0} \\
&=& 2\int_S \triangle_{g_0} u dV_{g_0} +2\int_S \frac{r}{2}e^{2u}-K_0 dV_{g_0}\\
&=& r( V(t)-V(0)).
\end{eqnarray*}
Here we have used Lemma \ref{lem_integration} and equation (\ref{cc}). Since $(V(t)-V(0))=0$ at $t=0$, it will remain so for $t\in [0,T]$.

\end{proof}

\section{Rotationally symmetric solitons}\label{five}
The purpose of this section is to show if $S$ is the Riemann sphere and $\beta$ consists of one or two conical singularities, then there exists (shrinking) Ricci Soliton metric in the conformal class $(S,\beta)$. As in the smooth case, a metric $g$ is called a gradient shrinking soliton if there exists some $f$ such that
\begin{equation}\label{eqn_defnsoliton}
	\frac{1}{2}R g_{ij}-\lambda g_{ij}=\nabla_{ij} f.
\end{equation}
\begin{rem}
	Here is $\lambda$ is some constant related to the volume. \textit{Assume} that Gauss-Bonnet theorem is true for this metric and $\int_S \triangle fdV=0$, then
	\begin{equation*}
		2\lambda Vol(S)=\int_S RdV=4\pi(\chi(S)+\sum \beta_i),
	\end{equation*}
	where $\chi(S)$ is the Euler number of $S$ and $\beta=\sum \beta_ip_i$. Therefore, one can take $\lambda=1$ by a scaling. From now on, we assume $\lambda=1$.
\end{rem}

\begin{prop}\label{prop_defnsoliton}
	Equation (\ref{eqn_defnsoliton}) is equivalent to

	1) $R-2=\triangle f$;

	2) The gradient field of $f$ is a conformal Killing field in the sense that the one parameter transformation group generated by $\nabla f$ is conformal.
\end{prop}

\begin{proof}
	Recall that a vector field $X$ is a conformal Killing field if and only if
	\begin{equation*}
		\nabla_i X_j+\nabla_j X_i=\frac{2}{n}\mbox{div} (X) g_{ij},
	\end{equation*}
	where $n$ is the dimensional of the manifold. In particular, 2) is equivalent to
	\begin{equation*}
		\nabla_{ij} f=\triangle f g_{ij}.
	\end{equation*}
	Here we used $n=2$. The proof is then trivial.
\end{proof}

From now on, we assume $S$ is the Riemann sphere and there is only one conical singularity $p$ of order $\beta$. One may take a conformal coordinate $(x,y)$ on $S$ such that $p$ is the infinity. The metric can be written as
\begin{equation*}
g=e^{2u}(dr^2+r^2 d\theta^2),
\end{equation*}
where $(r,\theta)$ is polar coordinates.

Let's assume that $g$ is rotationally symmetric such that $u$ is only a function of $r$. The condition that the order of the cone singularity is $\beta$ implies some asymptotic behavior of $u(r)$ when $r$ goes to infinity. To be precise, set $r=\tilde{r}^\sigma$ for some $\sigma<0$ to be determined. 
The metric is now
\begin{equation*}
g=e^{2u} \sigma^2 (\tilde{r}^{\sigma-1})^2 \left( d\tilde{r}^2 +\frac{1}{\sigma^2}\tilde{r}^2 d\theta^2 \right).
\end{equation*}
The metric has a cone singularity of order $\beta$ at $\tilde{r}=0$, which requires

1) $\sigma=-\frac{1}{\beta+1}$;

2) $e^{2u}\tilde{r}^{2\sigma-2}=e^{2u}r^{\frac{2\sigma-2}{\sigma}}$ has a finite limit when $r\to \infty$.
Hence,
\begin{equation}\label{eqn_asymp}
u(r)\sim -(\beta+2)\log r.
\end{equation}

Since we assume the soliton metric to be rotationally symmetric, it's natural that the conformal killing field in Proposition \ref{prop_defnsoliton} is also rotationally symmetric. It's not difficult to prove that on $\Real^2$, any rotationally symmetric conformal killing field is
\begin{equation*}
X=cr\pfrac{}{r},
\end{equation*}
where $c$ is some constant. By 2) in Proposition \ref{prop_defnsoliton},
\begin{equation*}
\nabla f= e^{-2u}\pfrac{f}{r}\pfrac{}{r}=cr\pfrac{}{r}.
\end{equation*}
Then
\begin{equation*}
\pfrac{f}{r}=cre^{2u}.
\end{equation*}
Hence,
\begin{equation}\label{eqn_laplacef}
\triangle f=e^{-2u}( f''+\frac{1}{r}f')=2c+2cru'.
\end{equation}

Recall that the scalar curvature 
\begin{equation}\label{eqn_R}
R=2e^{-2u}(-\triangle u)=-2e^{-2u}\left( u{''}+\frac{1}{r}u' \right).
\end{equation}

By the condition 1) in Proposition \ref{prop_defnsoliton} and equation (\ref{eqn_laplacef}) and (\ref{eqn_R}),
\begin{equation}
u''+(\frac{1}{r}+cre^{2u})u'+(1+c)e^{2u}=0.
\label{eqn_ODE}
\end{equation}

To find a soliton in the conformal class $(S,\beta)$, it suffices to find a solution to equation (\ref{eqn_ODE}) for \textit{some} $c$ such that $u'(0)=0$ and the asymptotic condition (\ref{eqn_asymp}) holds.

\begin{rem}\label{rem_c0}
	When $c=0$, equation (\ref{eqn_ODE}) has an explicit solution $u(r)=\log \frac{4}{4+r^2}$. This is nothing but the round sphere metric of radius one. The asymptotic behavior of $u(r)$ when $r\to \infty$ is $-2\log r$, which means the the cone singularity at infinity has order $\beta=0$. That is not really a singularity as can be seen from the round sphere.
\end{rem}

Equation (\ref{eqn_ODE}) is a Fuchsian ODE, or an ODE with a regular singular point. Thanks for Theorem 4.3 in Kichenassamy's book \cite{Kich}, the following local existence result is true.

\begin{thm}\label{thm_ode}
For each $c\in \Real$, there exists some $\varepsilon>0$ and $u(r)$ defined on $[0,\varepsilon]$ such that $u(r)$ is a solution to equation (\ref{eqn_ODE}) with initial value $u(0)=0$ and $u'(0)=0$. Moreover, the solution depends smoothly on $c$.
\end{thm}
\begin{proof}
Let $w=u'$. Equation (\ref{eqn_ODE}) can be written as a system
\begin{equation}
\label{eqn_fuchsian}
\begin{array}[]{lcl}
rw'+w &=& r(-cre^{2u}-(1+c)e^{2u}) \\
ru' &=& rw.
\end{array}
\end{equation}
This is the 'Fuchsian PDEs with analytic data' discussed by Kichenassamy in Section 4.2 of his book \cite{Kich}, if we take $c$ as a space variable. Theorem 4.3 in the same book implies that there exists a solution $u(r)$ defined on $[0,\varepsilon]$ for each $c$ and $u(r)$ depends smoothly on $c$. See Remark 4.2 in the same Section.
\end{proof}

We will then show for each $c$, the solution $u(r)$ is defined for $r>0$ and the asymptotic behavior of $u(r)$ as $r\to \infty$ is given by
\begin{equation*}
u(r)\sim -A_c \log r.
\end{equation*}
When $c$ varies in $(-1,+\infty)$, we will prove in the remaining part of this section that $A_c$ satisfies,

1) $A_c$ is a continuous function of $c$;

2) When $c$ goes to infinity, $A_c$ goes to $-1$;

3) When $c$ approaches $-1$ with $c>-1$, $A_c$ goes to $-\infty$.

Due to these properties of $A_c$, it's obvious that for each $\beta>-1$, there exists some $c$ such that
\begin{equation*}
A_c=-(\beta+2).
\end{equation*}
Hence, the solution to equation (\ref{eqn_ODE}) with this $c$ gives a Ricci soliton metric on $(S,\beta)$.

Let $u(r)$ be the solution given by Theorem \ref{thm_ode}. As long as the solution exists, define
\begin{equation}\label{eqn_AB}
A=ru',\quad B=re^{2u}.
\end{equation}
Equation (\ref{eqn_ODE}) is then equivalent to 
\begin{equation}
\left\{
\begin{array}[]{l}
A'=-B(cA+c+1) \\
B'=\frac{1}{r}B(2A+1).
\end{array}
\right.
\label{eqn_ode2}
\end{equation}

The initial condition for $u$ means $A(0)=0$ and $B(0)=0$. It turns out that the ODE system (\ref{eqn_ode2}) is easier to understand than equation (\ref{eqn_ODE}). Note that the ODE system (\ref{eqn_ode2}) is not of the type discussed in Theorem 4.3 of \cite{Kich}. The solution of the corresponding initial value problem may not be unique. In fact, it's not unique because the solution to equation (\ref{eqn_ODE}) with initial value $u(0)=C$ and $u'(0)=0$ gives rise to the different $A(r)$ and $B(r)$ with the same initial condition for different constant $C$. In the following proof, we will use the fact that $A$ and $B$ depends smoothly on the parameter $c$. Here the $A$ and $B$ are defined by equation (\ref{eqn_AB}) from $u$, which depends smoothly on $c$ as given in Theorem \ref{thm_ode}.

 Our next lemma proves that the the solution exists for all $r>0$ and when $r\to \infty$, $A(r)$ converges to some constant.

\begin{lem}\label{lem_limit}
	For all $c>-1$, the solution $A(r), B(r)$ to equation (\ref{eqn_ode2}) is defined on $[0,+\infty)$. Moreover, there exists some $A_c$ such that
\begin{equation*}
\lim_{r\to \infty} A(r)=A_c.
\end{equation*}
\end{lem}
\begin{proof}
	Claim: Under the assumption that $c>-1$, $A(r)$ is non-increasing as long as it exists.

	By the definition of $B$, we know $B\geq 0$. Due to the assumption $c>-1$, $c+1>0$. Hence, when $r$ is small and $A(r)$ small, $A'$ is nonpositive. $A$ will keep decreasing until $cA+c+1=0$, after which $A$ will be constant. Therefore, the claim is true.

	The proof of the lemma will be divided into three cases.

	In the case $c=0$, we know from Remark \ref{rem_c0} that the conclusion of the lemma is true and $A_0=-2$.

If $c>0$. $A$ will not stop decreasing before it hits
\begin{equation*}
-\frac{c+1}{c}.
\end{equation*}
Therefore, $A$ is bounded as long as it exists. This in turn implies that $B$ is bounded as long as it exists, because $B'=\frac{1}{r}B(2A+1)$ and both $A$ and $1/r$ are bounded for $r>\varepsilon$. Therefore, the solution exists for all $r>0$. $A(r)$ converges to some constant as $r$ goes to infinity.

Let's discuss the remaining case $0>c>-1$. If $A(r)>-2$ as long as it exists, then the same argument as above shows that the solution exists forever. Moreover the monotonicity of $A$ implies that $A(r)$ converges to some constant when $r\to \infty$. Hence we may assume at some time $r=r_0$. $A$ decreases to $-2$. For all $r>r_0$, we have
\begin{equation*}
B'< (-3)\frac{B}{r}.
\end{equation*}
ODE comparison gives
\begin{equation}\label{eqn_bbound}
B(r)<B(r_0)r_0^3 \frac{1}{r^3}
\end{equation}
for $r>r_0$.

On the other hand, since $A<-2$, we have some $\tilde{c}<0$ depends only on $c$, such that
\begin{equation*}
cA+c+1\leq \tilde{c}A.
\end{equation*}
Hence, for all $r>r_0$, we have
\begin{equation*}
A'> -\tilde{c} BA,
\end{equation*}
which implies
\begin{equation*}
(\log (-A))'<-\tilde{c} B.
\end{equation*}
So
\begin{equation*}
\log (-A(r))<\log (-A(r_0))+\int_{r_0}^r (-\tilde{c}) B(s)ds.
\end{equation*}
Equation (\ref{eqn_bbound}) then implies that $A(r)$ is bounded from below as long as it exists. The concludes the proof of the lemma.
\end{proof}

\begin{lem}\label{lem_kaka}
For $c>-1$, we have $A_c<-1$. Moreover, there exists some $r_1>0$ such that for $r>r_1$, $B(r)$ decreases and  $\lim_{r\to \infty}rB(r)=0$.
\end{lem}

\begin{proof}

We know $A(t)$ decreases. It suffices to show that $A(t)\geq -1$ for all $t$ is not possible. In that case,
\begin{equation*}
B'\geq - \frac{1}{r} B.
\end{equation*}
Hence, for some $r_0$ and $r>r_0$, we have
\begin{equation*}
B(t)\geq B(r_0)r_0\frac{1}{r}.
\end{equation*}
Since we assumed $A(t)\in [-1,0]$, there exists some constant $\kappa>0$ such that
\begin{equation*}
cA+c+1\geq \kappa.
\end{equation*}
In fact, one can take $\kappa=\frac{1+c}{2}$ if $0\geq c>-1$. If $c>0$, since $cA+c\geq 0$ due to our assumption, one may take $\kappa=1$. 

Therefore
\begin{equation*}
A'\leq -\kappa B\leq -\kappa B(r_0)r_0\frac{1}{r},
\end{equation*}
when $r>r_0$.
This is a contradiction to the assumption that $A\geq -1$ for all $r$.

Since $A_c$ will be smaller than $-1$, there exists some $r_1$ such that $2A(r_1)+1<A_c$. For $r>r_1$,
\begin{equation*}
B'=\frac{1}{r}B(2A+1)<\frac{A_c}{r}B.
\end{equation*}
Hence
\begin{equation*}
(\log B)'\leq \frac{A_c}{r}.
\end{equation*}
Integration over $r$ shows that there exists some constant $C$ depending on $r_1$ and $B(r_1)$ such that
\begin{equation*}
B(r)\leq Cr^{A_c}
\end{equation*}
when $r>r_1$. The second assertion of the lemma follows by observing that $A_c<-1$.
\end{proof}

The next lemma is the most important one.
\begin{lem}
As a function of $c$, $A_c$ is continuous for $c>-1$.
\end{lem}
\begin{proof}

For any $\varepsilon>0$, due to the previous lemma, there exists some $r_1>0$ such that

(1) $\abs{A(r_1)-A_c}\leq \varepsilon/10$;

(2) $r_1 B(r_1)\leq \eta\varepsilon$ for some $\eta>0$ to be determined later. $\eta$ will determined by nothing but $c$.

Assume without loss of generality that

(3) $A_c+\varepsilon<\frac{A_c+(-1)}{2}$.

For another parameter $\tilde{c}$, denote the solution to equation (\ref{eqn_ODE}) by $\tilde{u}$. Define $\tilde{A}$ and $\tilde{B}$ as in equation (\ref{eqn_AB}). Since $A$ and $B$ depend smoothly on $c$, if $\tilde{c}$ is sufficiently close to $c$, one may assume

(1') $\abs{\tilde{A}(r_1)-A_c}\leq \varepsilon/5$;

(2') $r_1\tilde{B}(r_1)\leq 2\eta \varepsilon$.

We will now estimate $\abs{A_{c'}-A_c}$. We will first prove some estimates under the condition that $\abs{\tilde{A}(r)-A_c}<\varepsilon$. These estimates together with (1') will ensure that $\abs{\tilde{A}(r)-A_c}<\varepsilon$ for all $r>r_1$, from which the lemma follows. 

If $\tilde{c}>0$, then $\tilde{c}\tilde{A}+\tilde{c}+1\leq \tilde{c}+1\leq c+2$. If $\tilde{c}\leq 0$, then $\tilde{c}\tilde{A}+\tilde{c}+1\leq 1+(-A_c+1)(\abs{c}+1)$, where we used the assumption that $\abs{\tilde{A}(r)-A_c}<\varepsilon$.
In either case, there exist some constant $\kappa$ depending only $c$ such that
\begin{equation*}
\tilde{c}\tilde{A}(r)+\tilde{c}+1\leq \kappa.
\end{equation*}
Hence
\begin{equation*}
\tilde{A}'\geq -\kappa \tilde{B}.
\end{equation*}

Integrate the above equation over $r$.
\begin{equation}\label{xx}
\tilde{A}(r)\geq \tilde{A}(r_1)-\kappa\int_{r_1}^r \tilde{B}(t)dt.
\end{equation}

Since $\abs{\tilde{A}(r_1)-A_c}<\varepsilon$, due to the monotonicity of $\tilde{A}$ and (3),
\begin{equation*}
\tilde{A}(r)<\frac{A_c-1}{2}.
\end{equation*}
Hence,
\begin{equation*}
\tilde{B}'(r)\leq \frac{A_c-1}{2r}\tilde{B}(r).
\end{equation*}
Integrate the above equation from $r_1$ to $r$,
\begin{equation}\label{yy}
\tilde{B}(r)\leq \tilde{B}(r_1)\frac{r_1^{(1-A_c)/2}}{r^{(1-A_c)/2}}.
\end{equation}
Plug equation (\ref{yy}) into equation (\ref{xx}),
\begin{equation*}
\tilde{A}(r)\geq \tilde{A}(r_1)-\kappa \frac{2}{-A_c-1}\tilde{B}(r_1)r_1.
\end{equation*}
Now we can choose $\eta$ so small that
\begin{equation*}
\kappa \frac{2}{-A_c-1}\eta<1/5.
\end{equation*}
By (1'),
\begin{equation*}
\tilde{A}(r)\geq A_c-\frac{3\varepsilon}{5}.
\end{equation*}
Hence our assumption that $\abs{\tilde{A}(r)-A_c}<\varepsilon$ is true for all $r>r_1$ and the lemma is proved.

\end{proof}

\begin{lem}
\begin{equation*}
\lim_{c\to \infty} A_c=-1
\end{equation*}
and
\begin{equation*}
\lim_{c\to -1}A_c=-\infty.
\end{equation*}
\end{lem}
\begin{proof}

For $c>0$, it's easy to see from the equation (\ref{eqn_AB}) that
\begin{equation*}
A(t)>- \frac{c+1}{c}.
\end{equation*}
The right hand side goes to $-1$ as $c$ goes to infinity. Combined with Lemma \ref{lem_kaka} and the monotonicity of $A(r)$, it follows that
\begin{equation*}
\lim_{c\to \infty} A_c=-1.
\end{equation*}

On the other hand, for $0>c>-1$,

For some $c$ very close to $-1$, by Lemma \ref{lem_kaka}, there exists $r_0>0$ such that $A(r_0)=-\frac{1}{2}$. Due to the monotonicity of $A$, for $r>r_0$, there exists some positive constant $\tilde{c}$, independent of $c$ such that
\begin{equation*}
cA+c+1>\tilde{c}(-A).
\end{equation*}
Equation (\ref{eqn_AB}) then implies
\begin{equation*}
(\log (-A))'>\tilde{c} B(r).
\end{equation*}

If the second assertion of this lemma is not true, i.e. there exists some $K>1$ such that for any $c$ and any $r>0$,
\begin{equation*}
A(r)>-K.
\end{equation*}

Then
\begin{equation*}
B'(r)\geq \frac{1}{r}B(-2K+1).
\end{equation*}
\begin{equation*}
B(r)\geq B(r_0)r_0^{2K-1}r^{-2K+1}.
\end{equation*}

\begin{equation*}
\log(-A(r))=\log(-A(r_0))+\tilde{c}\int_{r_0}^r B(t)dt.
\end{equation*}
The last integral is just a positive multiple of $B(r_0)r_0$. The proof is done if we can make $B(r_0)r_0$ as large as we want by choosing $c$ close to $-1$. Keep in mind that this $r_0$ depends on $c$ and it's determined by $A(r_0)=-\frac{1}{2}$.

To see this, we use the smooth dependence of ODE solution to the parameter again. Look at the equation where $c=-1$. Due to the uniqueness of the initial value problem to equation (\ref{eqn_ODE}), it's easy to check that if $c=-1$ then $u(r)\equiv 0$, hence $A\equiv 0$ and $B(r)=r$.

For any big $M$, we can choose $c$ so that
\begin{equation*}
A(M)<-1/2
\end{equation*}
and
\begin{equation*}
B(M)M>\frac{1}{2} M^2.
\end{equation*}

It's obvious from equation (\ref{eqn_AB}) that $B(r)$ is increasing as long as $A(r)\geq -\frac{1}{2}$. Due to Lemma \ref{lem_kaka} again, there is some $r_0>M$ such that $A(r_0)=-1/2$. Hence,
\begin{equation*}
B(r_0)r_0>\frac{1}{2}M^2.
\end{equation*}
This finishes the proof of this lemma.
\end{proof}

In summary, we have proved
\begin{thm}\label{thm_sub}
Let $S$ be the Riemann sphere and $\beta=\lambda p$ for some $p\in S$ and $\lambda>-1$. Then there exists some Ricci soliton metric in the conformal class $(S,\beta)$. 
\end{thm}

\begin{cor}
Let $S$ be the Riemann sphere and $\beta=\lambda_1p_1 + \lambda_2p_2$ for $p_1,p_2\in S$ and $\lambda_1,\lambda_2>-1$. Then there exists some Ricci soliton metric in the conformal class $(S,\beta)$.
\end{cor}

\begin{proof}
Set
\begin{equation*}
\lambda=\frac{\beta_1+1}{\beta_2+1}-1.
\end{equation*}
By Theorem \ref{thm_sub}, there exists some Ricci soliton metric $g$ in the conformal class $(S,\lambda p)$. Since the metric is rotationally symmetric, there are geodesics connecting $p$ and its antipodal point in every direction. Cut along two such geodesics such that their angle at $p$ is $\frac{2\pi (\beta_1+1)}{N}$ for some large $N$ to get a strip. Take $N$ copies of the strip and glue them side by side to get a rotationally symmetric metric with two conical singularities. By the definition of $\lambda$, the two conical singularities have the required cone angle.
\end{proof}

\begin{rem}
It's interesting to compare the Ricci soliton metric constructed above with the HCMU metric in \cite{HCMU}. The HCMU metric with two conical singularities, the football metric in \cite{HCMU}, may or may not have positive Gaussian curvature, depending on the ratio of the two cone angles. For the Ricci soliton metric constructed above,
\begin{eqnarray*}
K&=& e^{-2u}(-\triangle_0 u) \\
&=& -e^{2u}(u''+\frac{1}{r}u') \\
&=& -e^{2u}\frac{1}{r}A'(r).
\end{eqnarray*}
Recall that $A(r)=ru'(r)$. The monotonicity of $A$ implies $A'(r)<0$ for any $r$. Therefore, regardless of the cone angles, the Gaussian curvature of the Ricci soliton metric constructed above is positive. The author think it's reasonable to believe that this soliton metric is the least pinched metric on a tear drop or a football.
\end{rem}

\end{document}